\documentclass[11pt]{article}

	\usepackage{amsmath, amssymb, amscd, latexsym, graphicx}
\newcommand{\arrow}{\rightarrow}

\newcommand{\bb}{\mathbb}
\newcommand{\cx}{{\bb C}}

\newcommand{\half}{{\bb H}}
\newcommand{\integers}{{\bb Z}}

\newcommand{\reals}{{\bb R}}
\newcommand{\proj}{{\bb P}}

\newcommand{\makefig}[3]{
	\begin{figure}[htbp]
        \refstepcounter{figure}
	\label{#2}
        \begin{center}
		~#3~\\
		\medskip
                {\sf Figure \thefigure.  #1}
        \end{center}
	\end{figure}
}

\thicklines
\renewcommand{\bold}[1]{\smallskip \noindent {\bf \boldmath #1 }\nopagebreak[4]}

\newcommand{\qed}{\nopagebreak[4]\hfill
\rule{2mm}{2.5mm} \bigskip \pagebreak[2]}

\newcommand{\Mat}[1]{\left(\begin{matrix} #1 \end{matrix}\right)}

\renewcommand{\tilde}{\widetilde}

\newcommand{\asyto}{\sim}
\newcommand{\bab}{\asymp}

\newcommand{\bdry}{\partial}

\newcommand{\brackets}[1]{\langle #1 \rangle}

\newcommand{\closure}{\overline}

\newcommand{\disjunion}{\sqcup}
\newcommand{\bigdisjunion}{\bigsqcup}

\newcommand{\isom}{\cong}

\newcommand{\mem}{\in}

\newcommand{\plusorminus}{\pm}

\newcommand{\AND}{\;\;\;\text{and}\;\;\;}

\newcommand{\st}{\: : \:}         

\newcommand{\Hbar}{{\overline{H}}}

\newcommand{\Pbar}{{\overline{P}}}

\newcommand{\ML}{{\cal M}{\cal L}}
\newcommand{\ml}{\ML}

\renewcommand{\mod}{\operatorname{mod}}

\newcommand{\Mod}{\operatorname{Mod}}

\newcommand{\sys}{\operatorname{sys}}

\newcommand{\tw}{\operatorname{tw}}

\newcommand{\vol}{\operatorname{vol}}
\newcommand{\Vol}{\operatorname{Vol}}

\newcommand{\zed}{\integers}

\newcommand{\cC}{{\cal C}}

\newcommand{\cF}{{\cal F}}

\newcommand{\cM}{{\cal M}}

\newcommand{\cS}{{\cal S}}
\newcommand{\cT}{{\cal T}}

\newtheorem{theorem}{Theorem}[section]
\newtheorem{prop}[theorem]{Proposition}
\newtheorem{lemma}[theorem]{Lemma}
\newtheorem{cor}[theorem]{Corollary}

\makeatletter
\def\cleardoublepage{\clearpage\if@twoside \ifodd\c@page\else
    \thispagestyle{plain}\hbox{}\newpage\if@twocolumn\hbox{}\newpage\fi\fi\fi}

\def\ps@headings{\let\@mkboth\markboth
  \def\@oddfoot{}%
  \def\@evenfoot{}%
  \def\@evenhead{\small \sc\thepage\hfil\leftmark}
  \def\@oddhead{\small \sc \rightmark\hfil\thepage}
  \def\chaptermark##1{{
    \edef\@tempa{\ifnum \c@secnumdepth >\m@ne \@chapapp\ \thechapter. \fi}%
    \expandafter \markboth \expandafter{\@tempa ##1}{}}}%
  \def\schaptermark##1{\markboth {##1}{##1}}%
  \def\sectionmark##1{{
    \edef\@tempa{\ifnum \c@secnumdepth >\z@ \thesection. \fi}%
    \expandafter \markright \expandafter{\@tempa ##1}}}}

\def\thebibliography#1{\section*{References\@mkboth
 {References}{References}}\list
 {[\arabic{enumi}]}{\settowidth\labelwidth{[#1]}\leftmargin\labelwidth
 \advance\leftmargin\labelsep
 \usecounter{enumi}}
 \def\newblock{\hskip .11em plus .33em minus .07em}
 \sloppy\clubpenalty4000\widowpenalty4000
 \sfcode`\.=1000\relax}

\newif\if@restonecol
\def\theindex{\@restonecoltrue\if@twocolumn\@restonecolfalse\fi
\columnseprule \z@
\columnsep 35pt\twocolumn[\@makeschapterhead{Index}]
 \@mkboth{Index}{Index}\thispagestyle{plain}\parindent\z@
 \parskip\z@ plus .3pt\relax\let\item\@idxitem}
\def\@idxitem{\par\hangindent 40pt}

\def\endtheindex{\if@restonecol\onecolumn\else\clearpage\fi}

\def\footnoterule{\kern-3\p@ 
 \hrule width .4\columnwidth 
 \kern 2.6\p@} 
\@addtoreset{footnote}{chapter} 
\long\def\@makefntext#1{\parindent 1em\noindent 
 \hbox to 1.8em{\hss$^{\@thefnmark}$}#1}

\@addtoreset{equation}{section}

\renewcommand{\l@section}{\@dottedtocline{0}{1.5em}{2.3em}}
\renewcommand{\l@subsection}{\@dottedtocline{1}{3.8em}{3.2em}}
\renewcommand{\l@subsubsection}{\@dottedtocline{2}{7.0em}{4.1em}}

\makeatother

\newcommand{\dVol}{dV}
\newcommand{\Log}{\operatorname{Log}}

\newcommand{\alphat}{\tilde{\alpha}}
\newcommand{\gammat}{\tilde{\gamma}}

\input{xy}
\xyoption{all}

\title{ \vspace{-1in}
	\medskip
	{\bf  Teichm\"uller theory via random simple closed curves}
\vspace{.2in}}

\author{Curtis T. McMullen and Tina Torkaman}

\date{14 November 2025}

\begin{document}

\maketitle

\begin{abstract}
We show the map $\sigma : \cT_g \arrow \cC_g$
sending a compact hyperbolic surface $X$
to a random {\em simple} closed geodesic on $X$
determines a proper embedding of Teichm\"uller space
into the space of geodesic currents.
The proof depends on a formula for the
intersection number $i(C,C')$ 
of a pair of multicurves, expressed in terms of 
Dehn coordinates on $\ML_g(\zed)$.
\end{abstract}

\tableofcontents

\vfill \footnoterule \smallskip
{\footnotesize \noindent
        Research supported in part by the NSF.
	Typeset \today.
}

\thispagestyle{empty}
\setcounter{page}{0}

\newpage

\section{Introduction}
\label{sec:intro}

Let $\cT_g$ and $\cC_g$ denote the Teichm\"uller space of marked Riemann surfaces and the space
of geodesic currents attached to an oriented topological surface  $\Sigma_g$ of genus $g$.
We recall that $\cC_g$ is a {\em cone}:  it is closed under addition, and under multiplication by positive real numbers.
We let $i(C,C')$ denote the (geometric) intersection number between a pair of currents,
and $\ell_X(C)$ the hyperbolic length of $C$ on $X \mem \cT_g$.

The purpose of this note is to describe a new embedding
\begin{displaymath}
	\sigma : \cT_g \arrow \cC_g,
\end{displaymath}
based on random {\em simple} closed geodesics.  The map $\sigma$ can be compared to the traditional
embedding
\begin{displaymath}
	\lambda : \cT_g \arrow \cC_g,
\end{displaymath}
which is based on random {\em closed} geodesics.
We will see that $\sigma$ and $\lambda$ determine similar, but distinct, compactifications of $\cT_g$ 
by $\proj \ml_g$.  Along the way we obtain concrete estimates for lengths, intersection numbers
and volumes in Dehn coordinates on $\ML_g(\zed)$, of interest in their own right.

\bold{Random simple curves.}
The currents with $i(C,C) = 0$ form the space of {\em measured laminations} $\ML_g \subset \cC_g$.
Its integral points $\ML_g(\zed)$ are {\em multicurves}:  systems of disjoint simple closed curves with positive, integral weights.

We let $S_X$ denote the current describing a {\em random} simple closed curve on $X$.
It can be defined in two equivalent ways; first, as a limit
\begin{displaymath}
	S_X = \lim_{T \arrow \infty}  \frac{\sum S_X(T) }{T \,|S_X(T)|},
	\;\;\text{where}\;\;\;
	S_X(T) = \{C \mem \ML_g(\zed) \st \ell_X(C) < T\};
\end{displaymath}
and second, as an integral 
\begin{displaymath}
	S_X = \int_{B_X} C \,\dVol(C),
	\;\;\text{where}\;\;\;
	B_X = \{C \mem \ML_g(\reals) \st \ell_X(C) \le 1\} .
\end{displaymath}
Here $|S_X(T)| \bab T^{6g-6}$ is the number of elements in the 
finite set $S_X(T)$, and $\dVol$ is the natural (Thurston) measure on $\ML_g(\reals)$,
normalized so that $\ML_g(\zed)$ has density one (see equation (\ref{eq:Vol}) below).

Our main result (\S\ref{sec:main}) is:
\begin{theorem}
\label{thm:main}
The map $\sigma : \cT_g \arrow \cC_g$ defined by
\begin{displaymath}
	\sigma(X) = S_X/i(S_X,S_X)^{1/2}
\end{displaymath}
gives a proper embedding of Teichm\"uller space into the space of geodesic currents $C$ with $i(C,C)=1$.
\end{theorem}

\bold{Random closed geodesics.}
For comparison we recall that every $X \mem \cT_g$ also determines a
`random closed geodesic' $L_X \mem \cC_g$, its {\em Liouville current}, characterized by 
\begin{displaymath}
	i(C,L_X)  = \ell_X(C) 
\end{displaymath}
for all currents $C$; and the map $\lambda(X) = L_X$ (introduced in \cite{Bonahon:currents})
also gives a proper embedding
\begin{displaymath}
	\lambda : \cT_g \arrow \cC_g.
\end{displaymath}
This time the image lies on the locus $i(C,C) = \pi^2 (2g-2)$.
Certain features of Teichm\"uller space can be recovered from this embedding; for example,
the Weil--Petersson metric 
can be expressed in terms of the Hessian of $i(L_X,L_Y)$ along the diagonal \cite{Wolpert:Thurston:metric}.

\bold{The boundary of $\cT_g$.}
The closures of the images of $\lambda$ and $\sigma$ in $\proj \cC_g$ yield 
compactifications of $\cT_g$ by $\proj \ML_g$.  
In \S\ref{sec:pml} we compare these spaces, and show:

\begin{theorem}
\label{thm:pmlyes}
If $[C] \mem \proj \ML_g$ is a filling, uniquely ergodic lamination, and $\lambda(X_n) \arrow [C]$,
then $\sigma(X_n) \arrow [C]$.
\end{theorem}

\begin{theorem}
\label{thm:pmlno}
On the other hand, if $C = \sum_1^{3g-3} C_i$ is a filling multicurve, and $a_i > 0$ are real numbers,
then there exists a sequence $X_n \mem \cT_g$ such that
\begin{displaymath}
	\lambda(X_n) \arrow \left[\sum C_i \right]
	\;\;\;while\;\;\;
	\sigma(X_n) \arrow \left[\sum a_i C_i \right] .
\end{displaymath}
\end{theorem}
In other words, $\sigma(X_n)$ can converge to any point in the open simplex spanned by the $[C_i]$.
(Diagonalizing, one can deduce the same for the closed simplex.)

Thus $\sigma$ determines a compactification of $\cT_g$ that is
similar to, but distinct from, Thurston's original compactification
(treated in \cite{FLP}, \cite{Otal:book:fibered} and \cite{Bonahon:currents}).

\bold{Intersection numbers in Dehn coordinates.}
We now turn to an overview of the proof of Theorem \ref{thm:main}.

To study the current $S_X$, it is useful to introduce {\em Dehn coordinates}
$(m_i,t_i)_1^{3g-3}$ on the space of multicurves $\ML_g(\zed)$.
These coordinates record how simple curves intersect and twist relative to an
enhanced pants decomposition of $\Sigma_g$ (\S\ref{sec:Dehn}).  

One of the key ingredients in the proof is the following new topological result, of independent interest (\S\ref{sec:i}):
\begin{theorem}
\label{thm:i}
Let $C,C' \mem \ML_g(\zed)$ be multicurves with Dehn coordinates $(m_i,t_i)$ and $(m_i',t_i')$ respectively.
Then 
\begin{displaymath}
	i(C,C') = \sum_{i=1}^{3g-3} \left| \det \Mat{m_i & t_i \\m_i' & t_i'} \right|  + E,
\end{displaymath}
where $|E| \le N \sum_{i,j} m_i m_j'$ for some universal constant $N$.
\end{theorem}
The proof will show one can take $N=23$.  The statement extends without change to currents in $\ML_g(\reals)$.

\bold{Notation and conventions.}
Throughout the sequel we adopt the following conventions:
\begin{displaymath}
	\Log(x) = \max (1, \log x);
\end{displaymath}
the notation $A = O(B)$, or $A \ll B$, means $A < CB$ for some unspecified constant $C$, depending only on the genus $g$; and we write $A \bab B$ if $A \ll B \ll A$.

Let $\cM_g = \cT_g/\Mod_g$ denote the moduli space of Riemann surfaces of genus $g$.
The {\em systole} of $X \mem \cT_g$ is defined by:
\begin{displaymath}
	\sys(X) = \inf \{\ell_X(\alpha) \st \alpha \subset X \;\text{is a closed geodesic} \} .
\end{displaymath}
It is well--known that $[X_n] \arrow \infty$ in $\cM_g$ if and only if $\sys(X_n) \arrow 0$. 

\bold{Behavior of random simple curves.}
We will establish a sharper version of Theorem \ref{thm:i} in \S\ref{sec:i},
and review an estimate for $\ell_X(C)$ in Dehn coordinates in \S\ref{sec:L}.
When combined, these results yield:

\begin{cor}
\label{cor:X}
For all $X,Y \mem \cT_g$ we have:
\begin{eqnarray*}
	i(S_X,S_X) & \bab & \frac{\Vol(B_X)^2}{\sys(X) \Log(1/\sys(X))}, \;\text{and}\;\\[12pt]
	i(S_X,L_Y) & \gg  & C(Y) \Vol(B_X) / \sys(X) ,
\end{eqnarray*}
where $C(Y)>0$.
\end{cor}
For completeness we also recall that:
\begin{displaymath}
	\vol(B_X) \gg 1/\sys(X)^{1-\epsilon} ,
\end{displaymath}
an estimate we will review and sharpen in \S\ref{sec:L}; and note that
\begin{equation}
\label{eq:SXLX}
	i(S_X,L_X) = \frac{6g-6}{6g-5} \vol(B_X) 
\end{equation}
by a simple argument (integrate $\ell_X(C)$ over $B_X$).

\begin{cor}
The maps $\sigma_1(X) = S_X$ and $\sigma_2(X) = S_X/\Vol(B_X)$
also give proper embeddings of $\cT_g$ into $\cC_g$.
\end{cor}

\bold{Proof.}
In view of the bounds above, we have  
\begin{displaymath}
	i(S_X,S_X)^{1/2} \gg \Vol(B_X) \gg 1
\end{displaymath}
uniformly over $\cM_g$, and hence properness of $\sigma$ --- which is normalized
using the largest term above --- implies properness of the other two normalizations.
Continuity and injectivity are immediate from the same statements for $\sigma$;
in fact all 3 maps give the same embedding of $\cT_g$ into $\proj \cC_g$.
\qed

\bold{Properness of $\sigma$.}
To conclude our overview of the proof of Theorem \ref{thm:main},
we sketch the proof that $\sigma$ is a proper map.

Let $X_n \arrow \infty$ in $\cT_g$.  If $[X_n]$ stays in a compact subset of $\cM_g$, then we have $\sigma(X_n) \bab \lambda(X_n) \arrow \infty$
by properness of $\lambda$, and we are done.
On the other hand, if $X_n \arrow \infty$ in $\cM_g$, then $\sys(X_n) \arrow 0$,
and the estimates in Corollary \ref{cor:X} combine to give, for any fixed basepoint $Y \mem \cT_g$,
\begin{equation}
\label{eq:proper}
	i(\sigma(X_n),L_Y)^2  =
	\frac{i(S_X,L_Y)^2}{i(S_X,S_X)} 
	\gg \frac{C(Y)^2 \Log(1/\sys(X_n))
			}{
			\sys(X_n) 
		} \arrow \infty.
\end{equation}
Thus $\sigma(X_n) \arrow \infty$ in this case as well.  

For more details, see \S\ref{sec:main}.

\bold{Questions.}
\begin{enumerate}
	\item
How smooth is $i(S_X,S_Y)$ as a function on $\cT_g \times \cT_g$?
	\item
Can one define a metric on $\cT_g$ adapted to the embedding $\sigma : \cT_g \arrow \cC_g$?
	\item
For each $X \mem \cT_g$ we have a natural bijection $\cF : Q(X) \arrow \ML_g$,
sending a holomorphic quadratic differential $q$ to its horizontal measured foliation, $\cF(q)$.
Thus any probability measure $\mu$ on $Q(X)$ gives a current $S_X(\mu) = \int \cF(q) \cdot \mu$.

How do these currents behave when, for example, $\mu$ is normalized volume measure
on the unit ball $B \subset Q(X)$ for the Teichm\"uller metric?  For the Weil--Petersson metric?
\end{enumerate}

\bold{Notes and references.}
Our treatment of Dehn coordinates uses the theory of fractional Dehn twists,
presented in \cite{Luo:Stong:DT}.  The fact that maximally stretched geodesics are simple 
\cite{Thurston:stretch} also influenced the present discussion.
Another compactification of $\cT_g$, via extremal lengths, is discussed
in \cite{Gardiner:Masur:MF}.

\section{Winding numbers}
\label{sec:wind}

Let $X$ be a hyperbolic surface, and let $L(\delta)$ denote the 
hyperbolic length of a geodesic on $X$ (which may be an immersed
or embedded arc or loop).

In this section we will define the winding number
$w(\gamma,\alpha)$ for a geodesic arc $\gamma$ crossing a simple
closed geodesic $\alpha$ exactly once, and show:

\begin{theorem}
\label{thm:winding}
For any pair of arcs $\gamma,\gamma'$ joining the
two boundary components of a collar neighborhood $A = B(\alpha,r)$,
we have:
\begin{equation}
\label{eq:Lg}
	L(\gamma) = 2r + |w(\gamma,\alpha)| L(\alpha) + O(1),
\end{equation}
and
\begin{equation}
\label{eq:1}
	|\gamma \cap \gamma'| = |w(\gamma,\alpha) - w(\gamma',\alpha)| + E
\end{equation}
where $|E| \le 1$.
\end{theorem}
We will also show:

\begin{prop}
\label{prop:perp}
When $r=d(\alpha,\bdry A)$ is large, $\gamma$ is nearly perpendicular to $\bdry A$.
\end{prop}

\makefig{A pair of geodesics crossing in $S$.}{fig:crossing}{
\includegraphics[height=1.8in]{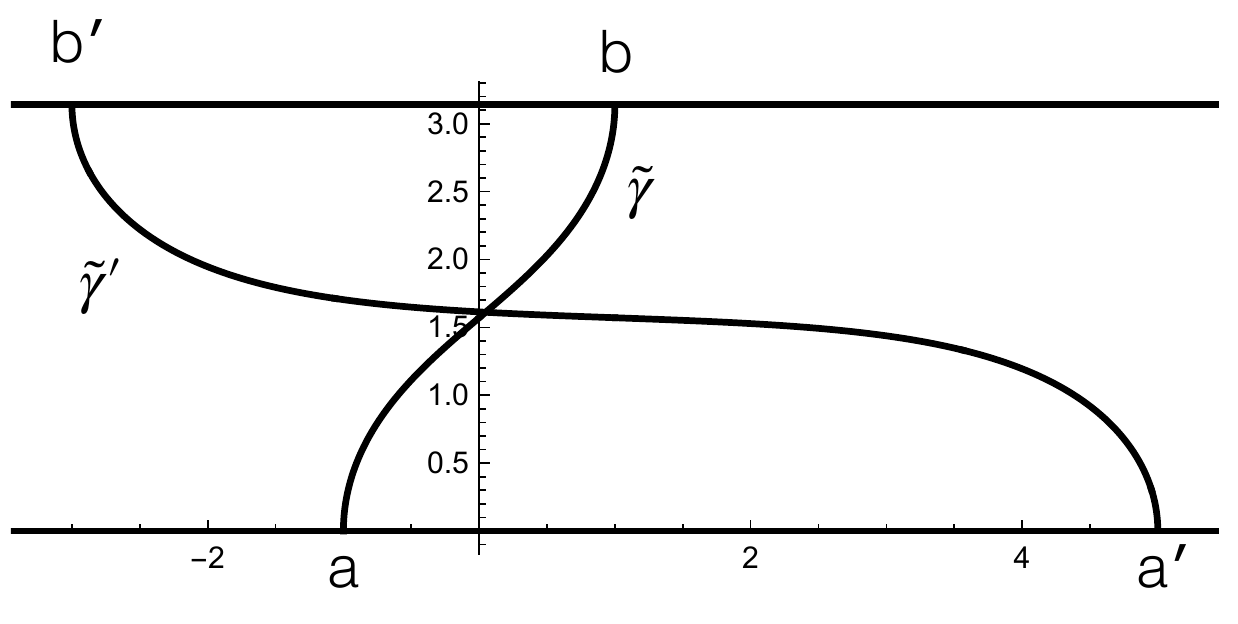}}

\makefig{A geodesic $\gamma \subset A$ with winding number $7/4$.}{fig:annulus}{
\includegraphics[height=2.2in]{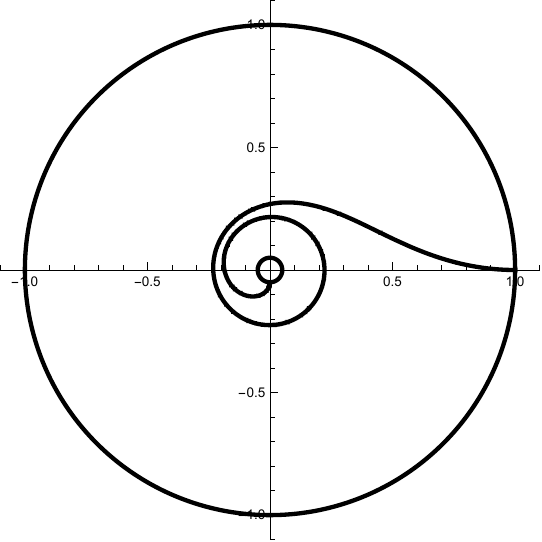}}

\bold{Geometric winding numbers.}
A {\em collar} about a simple closed geodesic $\alpha \subset X$ is an embedded annulus
of the form 
\begin{displaymath}
	A = B(\alpha,r) = \{x \mem X \st d(x,\alpha) < r\}.
\end{displaymath}
We allow $r=\infty$, in which case $A=X \isom S/\zed$, where 
\begin{equation}
\label{eq:S}
	S = \{z = x+iy \st 0<y<\pi\},\;\;\zed = \brackets{g},\AND g(z) = z+L.
\end{equation}
The geodesic $y=\pi/2$ in $S$ covers $\alpha$.

In this case, any complete geodesic $\gamma \subset A$ that crosses $\alpha$
is {\em embedded}.  Its lift to a geodesic $\gammat \subset S$
joins a pair of points of the form $a$ and $b+\pi i$ in $\bdry S$, and we define its
{\em geometric winding number} by
\begin{displaymath}
	w(\gamma,\alpha) = (b-a)/L .
\end{displaymath}
This real number is independent of the choice of lift and does not depend on
a choice of orientation for $\alpha$ or $\gamma$.
See Figure \ref{fig:crossing} (where we have used the map $\log : \half \arrow S$
to draw geodesics in $S$).

For an alternative definition, let $\pi : A \arrow \alpha$ denote
the nearest--point projection, whose fibers are geodesics orthogonal to $\alpha$.
Then $w(\gamma,\alpha)$ is the ratio between the immersed length of $\pi(\gamma) \subset \alpha$
and the length of $\alpha$ itself.  For an example with fractional 
winding number 
in the Poincar\'e metric on the annulus $A = \{z \st 1/20 < |z| < 1\} \subset \cx$, see Figure \ref{fig:annulus}.

\bold{General winding numbers.}
This second definition also makes sense for {\em any} geodesic segment $\gamma$ crossing $\alpha$ in a unique point $p$
on any hyperbolic surface $X$.  Indeed the orthogonal projection $\pi : \gamma \arrow \alpha$,
normalized so that $\pi(p)=p$, can be defined quite generally (e.g. by passing to the covering space $A = \half/\pi_1(\alpha)$ of $X$),
and then the winding number is given by:
\begin{displaymath}
	w(\gamma,\alpha) = \plusorminus \frac{L(\pi(\gamma))}{L(\alpha)}
\end{displaymath}
with the same sign convention as above.

\bold{Collars and crossings.}
In particular, for a finite collar of the form $A = B(\alpha,r) \subset X$ with $r<\infty$,
the winding number $w(\gamma,\alpha)$ is defined for any geodesic arc $\gamma \subset A$
that joins the two components of $\bdry A$.  

\bold{Proof of Theorem \ref{thm:winding}.}
The geodesic $\alpha$ cuts $\gamma$ into two segments of equal length $c$.
When lifted to the universal cover, each segment forms the hypotenuse of a right triangle
with sides of length $a=r$ and $b=w(\gamma)/2$.  It is well--known 
that the side lengths
of such a triangle satisfy $c = a+b + O(1)$ \cite[Theorem 2.2.2(i)]{Buser:book:spectra},
and (\ref{eq:Lg}) follows.

For the proof of equation (\ref{eq:1}), consider first the case where $r=\infty$ and $A=X$.
Write $A = S/\zed$ as in equation (\ref{eq:S}), and let $\gammat$ and $\gammat'$ be lifts of $\gamma$ and $\gamma'$ to $S$,
joining points $(a,b+\pi i)$ and $(a',b'+\pi i)$ on $\bdry S$ as in Figure \ref{fig:crossing}.
We have $|\gammat \cap \gammat'| = 1$ whenever $(a,b)$ and $(a',b')$ are linked, that is, whenever
$(a-a')(b-b') < 0$; otherwise, the intersection number is zero.  Thus $I = |\gamma \cap \gamma'|$ is the
same as the number of $n \mem \zed$ such that $(n+a-a')(b-b')< 0$.

Now it is easy to verify that for any interval $I = (x,x')$, we have $|I \cap \zed| = |x-x'| + E$
where $|E| \le 1$.  Since $w(\gamma) = b-a$ and $w(\gammat') = b'-a'$,
equation (\ref{eq:1}) follows.
\qed

\bold{Proof of Proposition \ref{prop:perp}.}
We refer again to Figure \ref{fig:crossing}, 
which shows the lift of $\gamma$ to a geodesic
$\gammat$ in the covering space $S$ of $X$ defined by 
$\pi_1(\alpha) \subset \pi_1(X)$.
The geodesic $\alphat \subset S$ defined $y=\pi/2$
covers $\alpha$, while the lift of $\bdry A$ to $S$ closest to $a$ is the horizontal
line $L$ defined by $y=\exp(-s)$, where $s = r + O(1)$.

Note that the angle between $\gamma$ and $\bdry A$ is the same
as the angle between $L$ and $\gammat$.
As $b \arrow \infty$ this angle gets smaller, and
$\gammat$ converges to the unique geodesic $\delta$ 
starting at $a$ and asymptotic to $\alphat$.
Since $\delta$ is orthogonal to $\reals$, it is nearly orthogonal
to $L$ when $r$ is large, and hence $\gamma$ is nearly
orthogonal to $\bdry A$ as well.
\qed

\section{Short geodesics}

In this section we will show:

\begin{prop}
\label{prop:thick}
Let $A  \subset X \mem \cM_g$ be the thin part of a compact Riemann surface with $3g-3$ short geodesics $(\alpha_i)$.
Let $G$ and $G'$ be a pair of simple closed geodesics on $X$, and let $\gamma$ and $\gamma'$ be connected components
of $G-A$ and $G'-A$ respectively.  Then:
\begin{displaymath}
	L(\gamma) = O(1) \AND |\gamma \cap \gamma'| \le 2,
\end{displaymath}
provided $\max L(\alpha_i)$ is sufficiently small.
\end{prop}

\bold{Collars and the thin part.}
We say a closed geodesic $\alpha_0 \subset X$ is {\em short} if $L = L(\alpha) \le 1$.
All short geodesics are simple.  The {\em standard collar} about a short geodesic is given by
\begin{displaymath}
	A_0 = B(\alpha_0,r), \;\;\text{where $\cosh(r) = 1/L$} .
\end{displaymath}
Each component of $A_0$ has length $1$.  This choice of $r$ guarantees that $A_0$ is an embedded cylinder,
and that the standard collars about distinct short geodesics are disjoint.
For proofs of these basic statements, see \cite[Ch. 4.1]{Buser:book:spectra}.
(We have used the fact that $1 < 2 \sinh^{-1}(1)$ and $\cosh^{-1}(1/L) < \sinh^{-1}(1/\sinh(L/2))$ for $L \le 1$
to simplify the statements above.)

The {\em thin part} of a compact Riemann surface $X \mem \cM_g$ whose short geodesics are $\alpha_1,\ldots,\alpha_n$ is defined by
the union of standard collars
\begin{displaymath}
	A = \bigcup_1^n A_i = \bigcup_1^n B(\alpha_i,r_i),
\end{displaymath}
where $\cosh(r_i) = 1/L(\alpha_i)$.  The maximum possible value of $n$ is $3g-3$, which is achieved
exactly when the curves $\alpha_i$ cut $X$ into $2g-2$ pairs of pants.

\makefig{The thick part $H$ of an ideal triangle, and the two types of geodesics arcs in the thick part of a pair of pants.}{fig:thin}{
\includegraphics[height=2.0in]{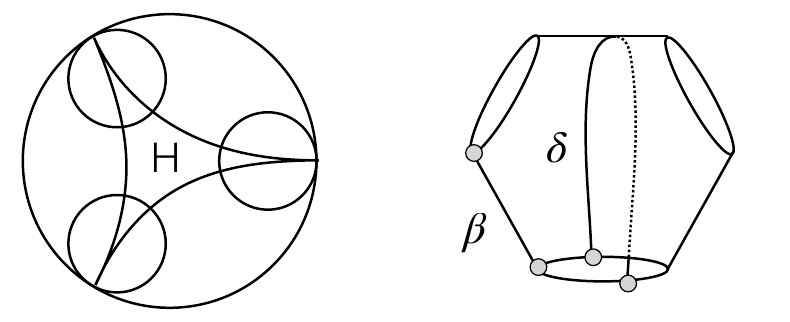}}

\bold{Proof of Proposition \ref{prop:thick}.}
Let $\alpha = \bigcup_1^{3g-3} \alpha_i$, and let $P$ denote the component of $X-\alpha$ containing $\gamma$.
Then $T = P-A$ is the thick part of the pair of pantsw $P$, and $\gamma \subset T$ is an arc joining two points on $\bdry T$.

Now assume that $\max L(\alpha_i)$ is very small.  Since the three boundary components of $P$
are contained in $\alpha$, they are very short, and hence $P$ itself is well--approximated by
a triply--punctured sphere.   It follows that $T$ is close to the double of the thick part $H$ of
an ideal triangle, shown at the left in Figure \ref{fig:thin}.
(Three of the boundary components of $H$ are horocycles of length $1/2$.)

Now by Proposition \ref{prop:perp}, $\gamma$ is nearly perpendicular to $\bdry T \subset \bdry A$.
Thus $\gamma$ is either close to the shortest geodesic $\beta$ joining a pair of distinct components
of $\bdry T$, or the shortest geodesic $\delta$ joining a single component to itself
(see Figure \ref{fig:thin}, right).  It follows readily that $L(\gamma) = O(1)$.
Two geodesics arcs of this type can meet in at most two points in $T$, and hence $|\gamma \cap \gamma'| \le 2$.
\qed

\section{Dehn coordinates}
\label{sec:Dehn}

In this section we review the definition of Dehn coordinates on $\ML_g(\zed)$.

We begin by defining a {\em hexagon decomposition} $(\alpha,\beta,H)$ of $\Sigma_g$.  This is an enhancement of a pants
decomposition defined by a multicurve $\alpha = \bigcup_1^{3g-3} \alpha_i$, allowing twists to be measured.
Parity conditions coming from the pants decomposition determine a lattice $L \subset (\zed^2)^{3g-3}$, 
with coordinates $(m_i,t_i)$, $i=1,\ldots,3g-3$.  A point $(m_i,t_i) \mem L$ is {\em positive} if for all $i$, either $m_i>0$ or $m_i=0$ and $t_i \ge 0$.
Letting $L_0 \subset L$ denote the set of positive parameters, we then define a multicurve 
$C = D(m_i,t_i)$ for each point $(m_i,t_i) \mem L_0$.  

These definitions provide the infrastructure that supports:

\begin{theorem}[Dehn]
\label{thm:Dehn}
The map $D : L_0 \arrow \ML_g(\zed)$ gives a bijection between positive
lattices points and isotopy classes of multicurves $C$ on $\Sigma_g$.
\end{theorem}
For the  proof, see the translation of Dehn's 1922 unpublished lecture notes, appearing in \cite{Dehn:sccs};
and \cite{FLP},  \cite{Harer:Penner:book} and \cite{Luo:Stong:DT}.
The fact that a pants decomposition must be enhanced to define Dehn coordinates is sometimes skirted in the literature.

\bold{Simple curves and (fractional) Dehn twists.}
Let $\cS_g$ denote the set of isotopy classes of unoriented, essential simple closed curves $\alpha \subset \Sigma_g$.
For convenience, we fix a hyperbolic metric on $\Sigma_g$ and identify $\alpha$ with its unique geodesic representative.
Since pairs of geodesics intersect efficiently, we have
\begin{displaymath}
	i(\alpha,\beta) = |\alpha \cap \beta|.
\end{displaymath}
for all $\alpha,\beta \mem \cS_g$.

\makefig{Twisting $\beta$ by $\alpha$.}{fig:astarb}{
\includegraphics[height=1.0in]{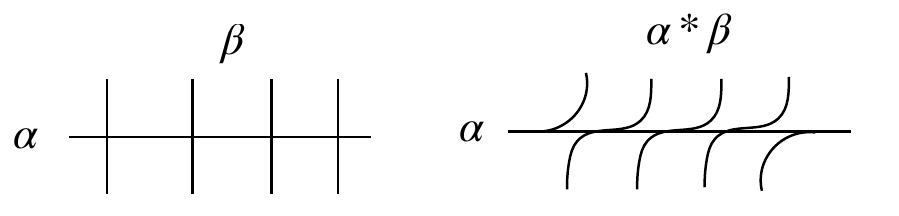}}

Each $\alpha \mem \cS_g$ determines a right Dehn twist $\tw_\alpha$ in the mapping--class group $\Mod_g$ of $\Sigma_g$.
Following \cite{Luo:Stong:DT}, we define a product on $\cS_g$ as follows:  provided $n = i(\alpha,\beta)>0$,
we define $\beta' = \alpha*\beta$ by modifying $\beta$ so that, as it crosses a collar neighborhood of $\alpha$,
its intersections with $\alpha$ shift one step to the right (see Figure \ref{fig:astarb}).
The result of iterating this operation $n$ times satisfies $(n \alpha)*\beta = \tw_\alpha(\beta)$,
so one can think of $\alpha*\beta$ as the result of the fractional Dehn twist $\tw_\alpha^{1/n}$.

If $i(\alpha,\beta)=0$, we set $\alpha*\beta = \beta$.

\bold{Hexagon decompositions.}
The definition of Dehn coordinates on $\ML_g$ depends on the choice of the topological data
$(\alpha,\beta,H)$, where:
\begin{enumerate}
	\item
$\alpha = \bigcup_1^{3g-3} \alpha_i$ is a maximal system 
of disjoint simple closed curves $\alpha_i \mem \cS_g$;
	\item
$\beta = \bigcup_1^n \beta_i$ is a second system of disjoint simple closed curves, with $i(\beta,\alpha_i) = 2$ for each $i$, such that 
$\Sigma_g-\beta$ is disconnected; and
	\item
$H$ is one of the two components of $\Sigma_g-\beta$.
The other component will be denoted $H^*$.
\end{enumerate}

\bold{Pairs of pants.}
The multicurve $\alpha = \bigcup_1^{3g-3} \alpha_i$ determines a {\em pants decomposition}
\begin{displaymath}
	\Sigma_g - \alpha = \bigdisjunion_{s=1}^{2g-2} P_s .
\end{displaymath}
Each $P_s$ is topologically a sphere with 3 disjoint disks removed.
Its boundary in $\Sigma_g$ has 3 or 2 boundary components.  In the first case, $\closure{P_s}$ is
an embedded pair of pants, while in the second, two cuffs are glued together.
In either case we write
\begin{equation}
\label{eq:Ps}
	\bdry P_s = \alpha_i \cup \alpha_j \cup \alpha_k,
\end{equation}
allowing for a repeated index when $\closure{P_s}$ has genus one.

\bold{Parity.}
The choice of $\alpha$ determines the lattice $L$:  it is the set of integers $(m_i,t_i)$ such that
\begin{displaymath}
	m_i + m_j + m_k = 0 \mod 2
\end{displaymath}
for each instance of equation (\ref{eq:Ps}).  It is easy to find the
{\em multiplicities} $(m_i)$ for a given multicurve $C$; they are determined by
$i(C,\alpha_i) = m_i$, and the parity condition reflects the fact that every arc in $C \cap P_s$    
has two endpoints.  The {\em twists} $t_i$ are more subtle are use the
curve $\beta$.

\makefig{Two hexagon decompositions $(\alpha,\beta,H)$ of $\Sigma_g$, $g=2$.  The curve $\beta$
has 3 components in the left example, and just 1 on the right.}{fig:hex}{
\includegraphics[height=1.7in]{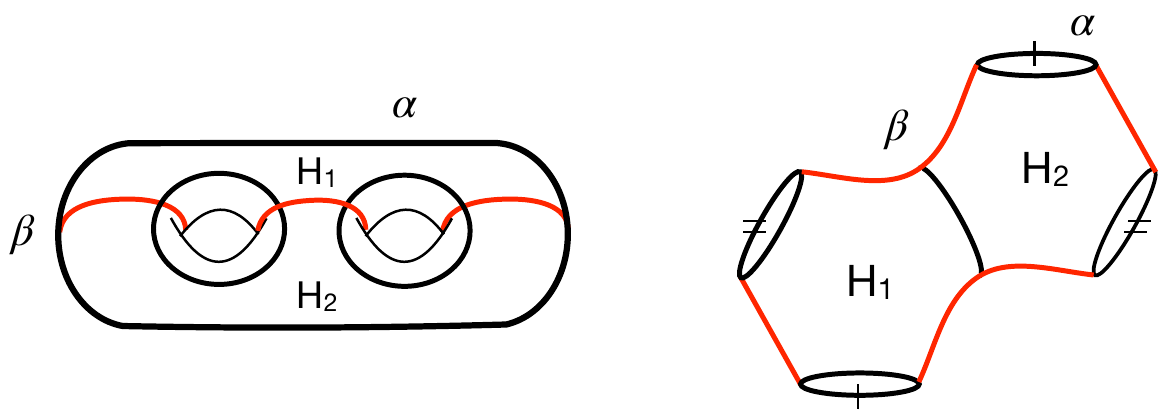}}

\bold{Hexagons.}
Each pair of pants is cut into two hexagons by $\beta$:  we have
\begin{displaymath}
	P_s - \beta = H_s \disjunion H_s^*,
\end{displaymath}
where $H_s = P_s \cap H$.  The sides of each hexagon alternate between $\alpha$ and $\beta$.
By construction, whenever $P_s$ and $P_t$ share $\alpha_i$ as a boundary component, $H_s$ and $H_t$ share the edge $\alpha_i \cap H$.

\bold{Examples.} 
Two examples of hexagon decompositions for $g=2$ are shown in Figure \ref{fig:hex}.
Note that in general, $\Sigma_g$ admits an orientation--reversing involution
$\iota$, fixing $\beta$ pointwise and satisfying $\iota(H) = H^*$.

\makefig{Arc systems $A(m)$ in $P$, with $(m_1,m_2,m_3)=(2,1,1)$ and $(4,1,1)$.}{fig:pants}{
\includegraphics[width=\textwidth]{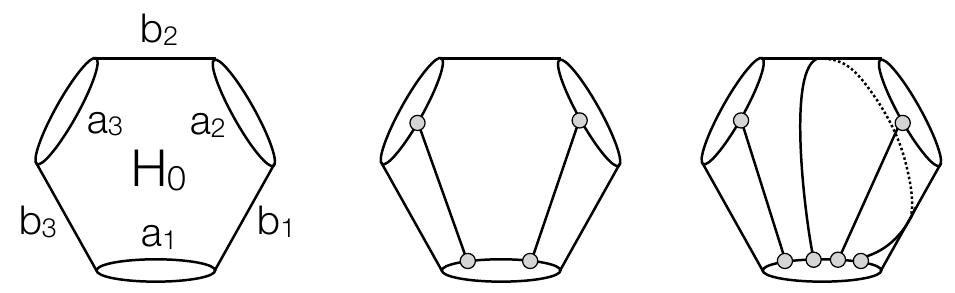}}

\bold{Arc systems.}
Let us now focus attention on a single, abstract pair of pants $P_0$,
decomposed into a pair of hexagons $H_0$ and $H_0^*$ meeting along
three edges, $b_1, b_2$ and $b_3$.  The remaining edges $a_1,a_2$ and $a_3$
of $H_0$ lie along the three components of $\bdry P_0$.
The indexing is chosen so the edges $(a_1,b_1,a_2,b_2,a_3,b_3)$ appear in counter--clockwise
order around $H_0$ (see Figure \ref{fig:pants}).

Let $m = (m_1, m_2, m_3)$ be a triple of integers such that $m_i \ge 0$ and $m_1+m_2+m_3$ is even.
We associate to $m$ a system of disjoint, properly embedded arcs 
$A(m) \subset \Pbar_0$
with the following properties:
\begin{itemize}
	\item
The arcs $A(m)$ have $m_i$ endpoints on $a_i \subset \bdry H_0$, and no other endpoints.
	\item
If an arc in $A(m)$ joins two distinct components of $\bdry P_0$, then it is contained in $\Hbar_0$.
	\item
If an arc in $A(m)$ joins $a_i$ to itself, then its interior crosses $\bdry H_0$ exactly
twice --- once at $b_i$, and once at $b_{i+1}$.  (Here the indexing is taken $\mod 3$.)
\end{itemize}
It is easy to show such an arc system exists.  
In fact, if the triangle inequality $m_i + m_j \ge m_k$ holds for every permutation $(ijk)$ of $(123)$,
then every arc joins a pair of {\em different} components of $\bdry P_0$, and we have $A(m) \subset \Hbar_0$.

On the other hand, if $m_i + m_j < m_k$, then every arc has at least one endpoint on $a_k$.
Any arc with {\em both} endpoints on $a_k$ passes once through $H_0^*$ before returning, and we have adopted the
{\em convention} (in the third bullet above) that in doing so it encircles $a_{k+1}$.  (Symmetry must be broken to insure
that the endpoints of $A(m)$ reside not just in $\bdry P_0$ --- a union of circles --- but in $\bdry P_0 \cap \bdry H_0$ ---
a union of arcs).

Using the fact that $H_0$ and $H_0^*$ are simply connected, one can verify:
\begin{quote}
	{\em The arc system $A(m)$ satisfying the conditions above is unique up to isotopy of
		the pair $(P_0,H_0)$.}
\end{quote}

\bold{Assembly and twists.}
It is now straightforward to define the multicurve $C = D(m_i,t_i)$.

First, suppose $C=D(m_i,0)$; that is, suppose $t_i=0$ for all $i$.
Choose a finite set $B \subset \alpha \cap \Hbar$
such that $|B \cap \alpha_i| = m_i$ for all $i$.
Next, for each pair of pants $P_s$ with boundary components $(\alpha_i,\alpha_j,\alpha_k)$, in counter--clockwise
order around $H_s$, construct the arc system $C_s = A(m_i,m_j,m_k) \subset P_s$ as above, 
adjusting by isotopy so that $\bdry C_s \subset B$.  Finally, we let
\begin{displaymath}
	C = D(m_i,0) = \bigcup C_s.
\end{displaymath}
This completes the definition of $D(m_i,t_i)$ when all $t_i=0$.

For general coordinates $(m_i,t_i) \mem L_0$, we simply define
\begin{displaymath}
	C = D(m_i,t_i) = \alpha_1^{t_1} * \alpha_2^{t_2} * \cdots * \alpha_{3g-3}^{t_{3g-3}} * D(m_i,0) + \sum_{m_j=0} t_j \alpha_j .
\end{displaymath}
The final sum accounts for the part of $C$ that does not cross $\alpha$, i.e.\ for curves parallel to $\bigcup \bdry P_j$.
This completes the definition of the mapping $D$ in Theorem \ref{thm:Dehn}.

\bold{Remarks.}
\begin{enumerate}
	\item
One can extend the definition of $D$ to the full lattice
$L$ by requiring that $D(\epsilon_i m_i, \epsilon_i t_i) = D(m_i,t_i)$, for any sequence $\epsilon_i \mem \{-1,1\}$.

	\item
It further extends to a continuous map
\begin{displaymath}
	D : (\reals^2)^{3g-3} \arrow \ML_g(\reals) , 
\end{displaymath}
generically of degree $2^{3g-3}$.

	\item
In Dehn coordinates, the natural (Thurston) measure on $\ML_g(\reals)$ is given locally by
\begin{equation}
\label{eq:dVol}
	\dVol = (1/V_g) \, dm_1 \, dt_1 \cdots dm_{3g-3} \,dt_{3g-3},
\end{equation}
where $V_g = [\zed^{6g-6}:L] = 2^{2g-2}$.  For any smoothly bounded open set $U \subset \ML_g$, we have
\begin{equation}
\label{eq:Vol}
	\Vol(U) = \lim_{T \arrow \infty} |T \cdot U \cap \ML_g(\zed)| / T^{6g-6} \cdot
\end{equation}
	
	\item
By smoothing $\beta \cup \alpha$, and potentially adding new branches to accommodate cases like the one at the right in Figure \ref{fig:pants},
one obtains a set of train tracks whose local charts covering the full space $\ML_g(\zed)$.
\end{enumerate}

\section{Intersection numbers}
\label{sec:i}

In this section we will prove and sharpen the estimate for $i(C,C')$ in Dehn coordinates provided by
Theorem \ref{thm:i}.  We will show:

\begin{theorem}
\label{thm:ii}
Let $(m_i,t_i)$ and $(m_i',t_i')$ denote the Dehn coordinates of $C,C' \mem \ML_g(\zed)$
with respect to a given a hexagon decomposition $(\alpha,\beta,H)$ of $\Sigma_g$.
We then have:
\begin{displaymath}
        i(C,C') = E + \sum_{i=1}^{3g-3} \max(0,|m_i t_i' - m_i' t_i| - 11 m_i m_i'),
\end{displaymath}
where 
\begin{displaymath}
	0 \le E \le 22 \sum m_i m_i' + \frac{1}{2}  \sum_{s=1}^{2g-2} (m_i + m_j + m_k)(m_i' + m_j' + m_k').
\end{displaymath}
Here $i,j,k$ and $s$ are related by $\bdry P_s = \alpha_i \cup \alpha_j \cup \alpha_k$.
\end{theorem}
We allow two of the indices $(i,j,k)$ to coincide when $\closure{P_s}$ is a torus with boundary.
The constants $11$ and $13$ are given only to show the bounds do not depend on $g$.

The idea of the proof is to realize $C$ and $C'$ as geodesics on a hyperbolic surface $X$
where $\alpha \perp \beta$ and $L(\alpha)$ is small.  This choice forces all of the
winding behavior into the thin part of $X$, as we will see below.

\bold{Fenchel--Nielsen coordinates from $(\alpha,\beta)$.}
Fix a hexagon decomposition $(\alpha,\beta,H)$ of $\Sigma_g$, with $\alpha = \bigcup_1^{3g-3} \alpha_i$.
This information determines, not just Dehn coordinates, but also a system of Fenchel--Nielsen coordinates $(\ell_i,\tau_i)$ on 
$\cT_g$ (see e.g.\ \cite[Chapter 4.6]{Thurston:book:TDGT}).  These coordinates are normalized so that 
\begin{displaymath}
	\ell_i = L(\alpha_i)
\end{displaymath}
for all $X \mem \cT_g$, and $\tau_i = 0$ for all $i$ if and only if $\alpha \perp \beta$ on $X$ --- meaning the
geodesic representatives of these curves are orthogonal at each of the $6g-6$ points where they meet.

\bold{Behavior of simple geodesics.}
We begin by describing the approximate geometry of (all) simple geodesics on $X$, at least when $\max \ell_i$ is small.

\begin{theorem}
\label{thm:G}
Suppose the Fenchel--Nielsen coordinates of $X \mem \cT_g$ satisfy $\tau_i = 0$ and $0 < \ell_i \ll 1$.  Let 
\begin{displaymath}
	A = \bigcup A_i = \bigcup B(\alpha_i,r_i)
\end{displaymath}
be the thin part of $X$, and let $G = G(m_i,t_i)$
be the geodesic representative of the multicurve $D(m_i,t_i)$,
Then the winding number of each component $\gamma$ of $G \cap A_i$ satisfies
\begin{equation}
\label{eq:5}
	|w(\gamma,\alpha_i) - t_i/m_i| \le 5,
\end{equation}
provided $m_i \neq 0$.
\end{theorem}

\bold{Proof.}
We simply repeat the construction of $C=D(m_i,t_i)$ in \S\ref{sec:Dehn}, this time using piecewise geodesic curves on $X$.

First suppose $t_i=0$ for all $i$.
For each $i$, choose $m_i$ distinct points on the arc $\alpha_i \cap H$.
Then for each pair of pants $P_s$, join these points by simple geodesics following the pattern 
shown in Figure \ref{fig:pants}.  The union of these arcs is a piecewise geodesic curve $G'$ 
representing $D(m_i,0) \mem \ML_g(\zed)$.

Note that each component $\gamma$ of $G' \cap A_i$ 
connects the two components of $\bdry A_i$ and meets $\beta$ at most once.
Thus it is nearly orthogonal to $\alpha_i$, and we have 
\begin{equation}
\label{eq:w1}
	|w(\gamma,\alpha_i)| \le 1.  
\end{equation}

Globally $G'$ is a finite union of polygonal loops.  Each loop potentially makes a small bend each time it crosses $\alpha$; otherwise it is straight,
and each maximal geodesic arc in $G'$ is long (since a collar neighborhood of $\alpha$ is wide).
It is a general principle in hyperbolic geometry that such a polygonal curve is $C^1$ close
to its geodesic representative $G$.  

In particular, the points where $G'$ crosses $\bdry A$ are very close to the points where $G$ does.
Since each component of $\bdry A$ has length one, the winding numbers of segments in $A$ change only by a small amount ($\ll 1$)
when we replace $G'$ with $G$.  Rounding this change up to $1$,
using equation (\ref{eq:w1}) and taking into account the fact that $A_i - \alpha_i$ has two components, we find that
\begin{equation}
\label{eq:G0}
	|w(\gamma,\alpha_i)| \le 3
\end{equation}
for each component $\gamma$ of $G \cap A_i$.  This completes the proof of Theorem \ref{thm:G} when all $t_i=0$.

Now let us incorporate twists.  Let $G_0 = G(m_i,0)$ be the geodesic
analyzed above.  Again following the discussion in \S\ref{sec:Dehn}, we 
construct a piecewise geodesic curve $G'$ by modifying $G_0$ inside $A_i$
as follows: first we replace $G_0$ with $\alpha_i^{t_i} * G_0$,
using the operation shown in Figure \ref{fig:astarb} in a small neighborhood of $\alpha_i$,
and then we straightening the resulting arcs to geodesics in $A_i$ while keeping their
endpoints on $\bdry A_i$ fixed.  Performing this operation for all $i$ leaves us with
a piecewise geodesic curve multicurve $G'$ such that:
\begin{enumerate}
	\item
$G'$ represents $D(m_i,t_i)$ in $\ML_g(\zed)$;
	\item
$G' = G_0$ outside $A$; and
	\item
$G'$ is a polygonal curve with small bends at the points $G' \cap \bdry A = G_0 \cap \bdry A$.
\end{enumerate}
The last point follows from Proposition \ref{prop:perp}.
Moreover, it is clear that if $t_i = nm_i$, then replacing $G_0$ with $G'$ changes the winding number
$w(\gamma,\alpha_i)$ by the integer $n$, because $\alpha_i^{nm_i} * G_0 = \tw_{\alpha_i}^{n}(G_0)$.  
 Since the winding number of each segment is monotone in $t_i$, it follows (using equation (\ref{eq:G0})) that each component $\gamma$ of $G' \cap A_i$ satisfies
\begin{displaymath}
	|w(\gamma) - t_i/m_i| \le 4.
\end{displaymath}
Using once more the fact that $G'$ is $C^1$ close to $G$, we find that the same inequality holds
with $G$ in place of $G'$, after relaxing the bound above to $5$.
\qed

\bold{Proof of Theorem \ref{thm:ii}.}
Choose $\ell_i=L$ small enough that both Proposition \ref{prop:thick} and Theorem \ref{thm:G} apply to $X$.
Let $G$ and $G'$ be the geodesic representatives of $C$ and $C'$ on $X$.
It suffices to treat the case $G$ and $G'$ and transverse to one another, as well as to $\alpha$, 
so that in particular we have
\begin{displaymath}
	i(C,C') = |G \cap G'|.
\end{displaymath}
To complete the proof, we will estimate the number of points of $G'\cap G$ in the thick and thin parts of $X$, then sum.

Let $T_s = P_s-A$ be the thick part of the pair of pants $P_s$. 
Let $\gamma$ and $\gamma'$ be components of $G \cap T_s$ and $G' \cap T_s$.
By Proposition \ref{prop:perp}, we have $|\gamma \cap \gamma'| \le 2$.
Summing over all pairs of pants, and using the fact
that the number of components of $G \cap T_s$ is $(m_i+m_j+m_k)/2$, we obtain the upper bound
\begin{equation}
\label{eq:GG1}
	|G \cap G' - A| \le (1/2) \sum_s (m_i+m_j+m_k)(m_i'+m_j'+m_k') ,
\end{equation}
where $\bdry P_s = \alpha_i \cup \alpha_j \cup \alpha_k$.

We now turn to the thin part $A = \bigcup A_i$.  Combining equations (\ref{eq:1}) and (\ref{eq:5}), we find that any pair of components
$\gamma$ and $\gamma'$ of $G \cap A_i$ and $G \cap A_i'$ satisfy:
\begin{displaymath}
	|\gamma \cap \gamma'| = \max(0,|t_i/m_i - t_i'/m_i'| + E),
\end{displaymath}
where $|E| \le 11$.  The number of such pairs is $m_i m_i'$; multiplying by this factor, summing over $i$, and using
the intermediate value theorem, this gives
\begin{equation}
\label{eq:GG2}
	|G \cap G' \cap A| = \sum_i \max(0, |m_i't_i - m_i  t_i'|+E m_i m_i'),
\end{equation}
with the same bound on $E$.  Adding equations (\ref{eq:GG1}) and (\ref{eq:GG2}), we obtain
the statement of the Theorem.
\qed

\bold{Remark.}
The strategy of using particular hyperbolic structures to prove a topological theorems on surfaces is
also used, for example, in \cite{Chas:McMullen:Phillips:sccs}.
For more estimates on lengths and intersection numbers, see \cite{Torkaman:interaction}.

\section{Lengths}
\label{sec:L}

In this section we give an estimate for hyperbolic length in Dehn coordinates, similar to that stated in \cite{Mirzakhani:thesis1} as Proposition 3.5.
We will then use this estimate to give a good qualitative picture of $S_X$ in Dehn coordinates.

\bold{Lengths.}
Fix a hexagon decomposition $(\alpha,\beta,H)$ for $\Sigma_g$, with $\alpha = \bigcup \alpha_i$, and let
$D(m_i,t_i)$ denote the corresponding Dehn coordinates on $\ML_g$.
We will show:

\begin{theorem}[Mirzakhani]
\label{thm:M}
For all $X \mem \cT_g$ there exists a hexagon decomposition
$(\alpha,\beta,H)$ of $\Sigma_g$ such that in Dehn coordinates,
\begin{equation}
\label{eq:M}
	\ell_X(D(m_i,t_i)) \bab \sum_1^{3g-3} 2 m_i \Log(1/\ell_i) + |t_i| \ell_i ,
\end{equation}
where $\ell_i = \ell_X(\alpha_i)$.  
\end{theorem}

The proof we give yields the following complement:
\begin{theorem}
\label{thm:Msharp}
Suppose $X \mem \cT_g$ has Fenchel--Nielsen coordinates $(\ell_i,\tau_i)$ with respect to $(\alpha,\beta,H)$,
and $\tau_i=0$ for all $i$.  Then the implicit constants in equation (\ref{eq:M}) tend to $1$ as $\max \ell_i \arrow 0$.
\end{theorem}

\bold{Adaptation.}
When the hexagon decomposition $(\alpha,\beta,H)$ is chosen so that Theorem \ref{thm:M} holds, we say the
corresponding Dehn coordinates are {\em adapted} to $X$.
As we will see in the proof,
one can always assume that $\alpha$ contains all short geodesics on $X$,
that $L(\alpha) = O(1)$, and that there is a definite angle between $\alpha$
and $\beta$ wherever they meet.
In Theorem \ref{thm:Msharp}, $\beta$ and $\alpha$ are orthgonal on $X$.

\bold{Approximating $B_X$ by a cube.}
Theorem \ref{thm:M} implies that $B_X$ is approximately a cube.
More precisely, if we let
\begin{equation}
\label{eq:CX}
        C_X = \prod_1^{3g-3} \left[0,\frac{1}{\Log(1/\ell_i)}\right] \times \left[\frac{-1}{\ell_i},\frac{1}{\ell_i}\right] \subset \ML_g,
\end{equation}
where $\ell_i = \ell_X(\alpha_i)$, then we have:

\begin{cor}
\label{cor:CX}
There exists a constant $k \bab 1$ such that for all $X \mem \cT_g$, we have
\begin{displaymath}
        (1/k) C_X \subset B_X \subset k C_X 
\end{displaymath}
in Dehn coordinates adapted to $X$.
\end{cor}

\begin{cor}
\label{cor:VolBX}
We have $\Vol(B_X) \bab \prod_1^{3g-3} (\ell_i \Log(1/\ell_i))^{-1}$.
In particular, $\Vol(B_X) \arrow \infty$ as $[X] \arrow \infty$ in $\cM_g$.
\end{cor}
This Corollary of Theorem \ref{thm:M} is also stated in \cite[Theorem 1.5]{Arana:Athreya:square}.

\bold{Proof of Theorems \ref{thm:M} and \ref{thm:Msharp}.}
Let $G \subset X$ denote the geometric representative of $C$, with hyperbolic length $L(G)$.
We may assume that each component of $G$ has multiplicity one, so that $L(G) = \ell_X(C)$.

First assume that $X$ has Fenchel--Nielsen coordinates $(\ell_i,\tau_i)$ with $\tau_i = 0$ and $\epsilon = \max \ell_i$ is very small --- in particular,
small enough that Proposition \ref{prop:thick} applies.
We will write $a \asyto b$ if $a/b \arrow 1$ as $\epsilon \arrow 0$.

Let $A = \bigcup A_i = \bigcup B(\alpha_i,r_i)$ be the thin part of $X$.  Then $\cosh(r_i) = 1/\ell_i$, so $r_i \asyto \Log(1/\ell_i)$.   
Combining Theorem \ref{thm:winding} on length with Theorem \ref{thm:G} on winding numbers, we obtain
\begin{displaymath}
	L(G \cap A) = \sum 2 m_i \Log(1/\ell_i) + |t_i| \ell_i + O(m_i) .
\end{displaymath}
Theorem \ref{thm:G} also shows the number of components of $G-A$ is $O\left(\sum m_i\right)$.
By Proposition  \ref{prop:thick}, each component of $G-A$ has length $O(1)$, and hence
\begin{displaymath}
	L(G-A) = O\left(\sum m_i\right) 
\end{displaymath}
as well.  Since $\Log(1/\ell_i) \arrow \infty$ as $\ell_i \arrow 0$, we can neglect terms of size $O(1)$ and hence 
\begin{displaymath}
	L(G) \asyto \sum 2 m_i \Log(1/\ell_i) + |t_i| \ell_i 
\end{displaymath}
as desired.  

This proves complete the proof of Theorem \ref{thm:Msharp}.   It shows there exists an $\epsilon_0>0$
such that Theorem \ref{thm:M} holds for $X$ provided its Fenchel--Nielsen coordinates satisfy $\tau_i = 0$ for all $i$ and $\max \ell_i < \epsilon_0$.

To handle the case of a general Riemann surface of genus $g$, first recall
that any $X \mem \cT_g$ admits a pants decomposition defined by $\alpha = \bigcup_1^{3g-3} \alpha_i$
with $\max L(\alpha_i) = O(1)$.  (See the discussion of Bers' constant in \cite[Ch. 5]{Buser:book:spectra}.
The implicit constant depends on $g$, as is permitted by our conventions from \S\ref{sec:intro}.)
We can then complete $\alpha$ to a hexagon decomposition $(\alpha,\beta,H)$ of $X$.

Let $(\ell_i,\tau_i)$ be Fenchel--Nielsen coordinates for $X$ compatible with $(\alpha,\beta,H)$.
Then $\ell_i=O(1)$, and we can modify $\beta$ by Dehn twists along components of $\alpha$ so that
$\tau_i = O(\ell_i)$.

Now let $X'\mem \cT_g$ be the unique Riemann surface with coordinates Fenchel--Nielsen coordinates $(\ell_i',\tau_i')$, where
$\ell_i' = \min(\epsilon_0,\ell_i)$ and $\tau_i'=0$.  Since $\ell_i' \le \epsilon_0$ for all $i$,
(\ref{eq:M}) holds for $X'$.  On the other hand, the Fenchel--Nielsen coordinates for $X$ and $X'$
are close enough that there exists a $K$--bilipschitz map $\phi : X \arrow X'$,
compatible with markings and with $K=O(1)$.
This mapping changes both sides in equation (\ref{eq:M}) by at most a bounded factor.
Since we have already proved Theorem \ref{thm:M} for $X'$, this shows that it also holds (with a bounded change of constants)
for general $X \mem \cT_g$.
\qed

\bold{Remark.}
The discussion in \cite[\S 3]{Mirzakhani:thesis1} omits some significant details regarding twists.

\section{Random simple closed curves}
\label{sec:X}

In this section we will use the estimates on $i(C,C')$ and $\ell_X(C)$ from \S\ref{sec:i} and \S\ref{sec:L} to prove
Corollary \ref{cor:X}, which states that for all $X, Y \mem \cT_g$ we have
\begin{eqnarray*}
        i(S_X,S_X) & \bab & \Vol(B_X)^2 (\sys(X) \Log(1/\sys(X)))^{-1} , \;\text{and}\;\\
        i(S_X,L_Y) & \gg  & C(Y) \Vol(B_X) / \sys(X) .
\end{eqnarray*}
For the proofs we study the interaction between simple closed curves
chosen at random in $B_X$.

\bold{Random pairs of simple closed curves.}
A {\em random variable} $V$ on a finite measure space $(M,\mu)$ is simply a measurable function $V : M \arrow \reals$;
its {\em expectation} $E(V)$ is given by $(1/\mu(M)) \int_M V\, d\mu$. 

Given $X \mem \cT_g$ we will be interested in currents $C$ chosen at random with respect to the natural volume element $\dVol$ on $B_X$
defined by equation (\ref{eq:dVol}).
In view of Corollary \ref{cor:CX}, $C_X$ is a good approximation
to $B_X$.  The advantage is that, in view of the
product structure on $C_X$ given by equation (\ref{eq:CX}), 
the coordinates $(m_i,t_i)$ of a random
current $C \mem C_X$ are {\em independent random variables},
each uniformly distributed in the corresponding interval:
\begin{equation}
\label{eq:ind}
	m_i \mem [0,1/\Log(1/\ell_i)]
	\AND
	t_i \mem [-1/\ell_i,1/\ell_i].
\end{equation}
Using this idea, the following two estimates lead to a proof of Corollary \ref{cor:X}.

\begin{theorem}
\label{thm:CC}
For $C,C' \mem B_X$ chosen at random, we have
\begin{displaymath}
	E(i(C,C')) \bab (\sys(X) \Log 1/\sys(X))^{-1} .
\end{displaymath}
\end{theorem}

\bold{Proof.}
By compactness it suffices to prove the result when $\sys(X)$ is small.  We will work in coordinates provided by a hexagon decomposition $(\alpha,\beta,H)$
adapted to $X$.   We can assume that $\alpha_1$ is the shortest closed geodesic on $X$, and hence
the Fenchel--Nielsen coordinates $(\ell_i,\tau_i)$ of $X$ satisfy $\ell_1=\sys(X)$.  
Moreover, by Corollary \ref{cor:CX}, we can replace $B_X$ with $C_X$ in the statement of the Theorem.
Since the coordinates of a pair of randomly chosen points $C,C' \mem C_X$
are independent random variables, chosen uniformly in
the intervals given by equation (\ref{eq:ind}),
we have
\begin{displaymath}
	E \left|\det \Mat{m_i & t_i \\ m_i' & t_i'} \right| \bab E(m_i)E(|t'_i|) \bab 
	\frac{1}{\ell_i \Log(1/\ell_i)} \cdot
\end{displaymath}
For $i=1$ this final quantity is $(\sys(X) \Log 1/\sys(X))^{-1})$; while for other values of $i$,
it is only smaller.  Using Theorem \ref{thm:i} we find that
\begin{displaymath}
	E(i(C,C')) \bab  (\sys(X) \Log 1/\sys(X))^{-1} + E,
\end{displaymath}
where $|E| \ll (\Log(1/\sys(X)))^{-2} $.
The error term $E$ is much smaller than the main term when $\sys(X)$ is small, and the Theorem follows.
\qed

\begin{theorem}
\label{thm:LY}
Fix $Y \mem \cT_g$, and let $C \mem B_X$ be chosen at random.  We then have
\begin{equation}
\label{eq:ell}
	E(i(C,L_Y)) \gg C(Y) / \sys(X) . 
\end{equation}
\end{theorem}

\bold{Proof.}
It suffices to prove this inequality when $\sys(X)$ is small.
Otherwise, we have $i(S_X,L_Y) \bab i(L_X,L_Y) \ge i(L_Y,L_Y) \bab 1$.
(The inequality $i(L_X,L_Y) \ge i(L_Y,L_Y)$ follows when $X$ is close to $Y$ by \cite[eq. (0.1)]{Wolpert:Thurston:metric},
and then for all $X$ and $Y$ by convexity of lengths along earthquake paths \cite[Theorem 1]{Kerckhoff:Nielsen}.)

By Corollary \ref{cor:CX}, it suffices to prove the statement with $B_X$ replaced by $C_X$.
Choose a system of simple closed geodesics $\gamma = \bigcup_1^n \gamma_i$
that fill the surface $Y$, in the sense that every component of $Y-\gamma$
is a topological disk.  
The current $L_Y$ is well--approximated by $\sum \gamma_i$; in particular
there is a constant $C(Y)>0$ such that
\begin{equation}
\label{eq:gamma}
	i(C,L_Y) \ge C(Y) \max_i i(\gamma_i,C) 
\end{equation}
for all $C \mem \cC_g$.

Let us examine the interaction of the curves $(\gamma_i)$ with
the loops $(\alpha_i)$ giving Dehn coordinates on $X$,
ordered so that $\ell_1 = \ell_X(\alpha_1) = \sys(X)$.
Since $\gamma$ fills $Y$, we have $i(\gamma_i,\alpha_1) \ge 1$ for some $i$; we may
assume that $i=1$.  Thus in the Dehn coordinates determined by $(\alpha,\beta,H)$,
we have $\gamma_1 = (m_i',t_i')$ with $m_i' \ge 1$.

Since every term in the sum appearing in Theorem \ref{thm:ii} is positive, for all $C = (m_i,t_i) \mem C_X$ we have
\begin{displaymath}
	i(C,\gamma_1) \ge |t_1 m_1' - t_1'm_1| - 11 m_1 m_1'.
\end{displaymath}

Using again the fact that coordinates $(m_1,t_1)$ of a random point
$C \mem C_X$ are independent random variables,
we find:
\begin{displaymath}
	E[|t_1 m_1' - t_1'm_1|] \gg  m_1'(1/\ell_1) ,
\end{displaymath}
while
\begin{displaymath}
	E(m_1 m_1') \le m_1' /\Log(1/\ell_1) .
\end{displaymath}
Once $\ell_1 = \sys(X)$ is small enough, the first expression
is the dominant one, and since $m_1' \ge 1$ we find
\begin{displaymath}
	E(i(C,\gamma_1)) \gg 1/\sys(X) .
\end{displaymath}
The Theorem now follows from equation (\ref{eq:gamma}).
\qed

\bold{Proof of Corollary \ref{cor:X}}.
The desired estimates follow from Theorems \ref{thm:CC} and \ref{thm:LY} above,
using the fact that for $C$ chosen randomly in $B_X$ and any $C' \mem \cC_g$,
we have
\begin{displaymath}
	E(i(C,C')) = i(S_X,C')/\Vol(B_X) .
\end{displaymath}
\qed

\section{Embedding of $\cT_g$}
\label{sec:main}

In this section we prove Theorem \ref{thm:main}, which states that the map
\begin{displaymath}
	\sigma : \cT_g \arrow \cC_g
\end{displaymath}
given by $\sigma(X) = S_X/i(S_X,S_X)^{1/2}$ is a proper embedding.

\begin{lemma}
\label{lem:K}
For any compact subset $K \subset \cM_g$, there exists a constant $C(K)>0$ such that
\begin{displaymath}
	i(\sigma(X),C) \ge C(K) \,i(\lambda(X),C) 
\end{displaymath}
for all $X \mem K$ and $C \mem \cC_g$.
\end{lemma}

\bold{Proof.}
Note that the function
\begin{displaymath}
	f(X,C) = \frac{i(\sigma(X),C)}{i(\lambda(X),C)}
	=
	\frac{i(\sigma(X),C)}{\ell_X(C)}
\end{displaymath}
is continuous $\cT_g \times \proj \cC_g$,
and its inf over $\cC_g$ descends to a continuous function of $[X]$ on $\cM_g$.
The result then follows from compactness of $K$.
\qed

\bold{Proof of Theorem \ref{thm:main}.}
Continuity of $\sigma$ follows from the continuity of $\ell : \cT_g \times \cC_g \arrow \reals$,
and injectivity follows from \cite[Theorem 1.4]{Erlandsson:Souto:distr} 
(In brief, from the ergodic decomposition of the current $S_X$ one can
recover $B_X$ up to a scalar factor, and this information determines $X$ uniquely
by \cite[Prop. 3.5]{Thurston:stretch}.)
It remains to prove that $\sigma$ is proper.

Fix a basepoint $Y \mem \cT_g$, and suppose $X_n \arrow \infty$ in $\cT_g$.
We treat two different cases.

First, suppose $[X_n]$ ranges in a compact subset $K \subset \cM_g$.
Then by Lemma \ref{lem:K} with $C=L_Y$, we have
\begin{displaymath}
	i(\sigma(X_n),L_Y) \ge C(K) i(\lambda(X_n),L_Y) \arrow \infty
\end{displaymath}
by properness of $\lambda$.  Thus $\sigma(X_n) \arrow \infty$ as desired.

Second, suppose $[X_n] \arrow \infty$ in $\cM_g$.  Then $\sys(X_n) \arrow 0$.
Squaring the second bound in Corollary \ref{cor:X} and dividing by the first, with $X=X_n$, we obtain
\begin{displaymath}
	i(\sigma(X_n),L_Y)^2 \gg \frac{\Log (1/\sys(X_n))}{\sys(X_n)} \arrow \infty,
\end{displaymath}
so $\sigma(X_n) \arrow \infty$ in this case as well (cf.\ equation 
(\ref{eq:proper}) from the Introduction).

Since we can always pass to a subsequence such that either $[X_n]$ is bounded in $\cM_g$ or $[X_n] \arrow \infty$,
the proof is complete.
\qed

\section{Compactification}
\label{sec:pml}

In this section we prove Theorems \ref{thm:pmlyes} and
\ref{thm:pmlno}, describing the behavior near a given
$[C] \mem \proj \ML_g$ of the two compactifications of
$\cT_g$ determined by $\lambda$ and $\sigma$.

\begin{lemma}
If $\lambda(X_n) \arrow [C] \mem \proj \ML_g$, and $\sigma(X_n) \arrow [C']$,
then $i(C,C') = 0$.
\end{lemma}

\bold{Proof.}
Fix a basepoint $Y \mem \cT_g$ and let $U \subset \cC_g$ be the compact set of currents
defined by $i(C,L_Y) = 1$.  We have a natural maps
\begin{displaymath}
	\cC_g-\{0\} \arrow \proj \cC_g \isom U
\end{displaymath}
whose composition is given by $C \mapsto C^u = C/i(C,L_Y)$.  Using the fact that $\sys(X) = O(1)$ on $\cM_g$,
equation (\ref{eq:SXLX}) and Corollary \ref{cor:X} 
give $i(S_X,L_X) \bab \vol(B_X)$
and $\vol(B_X) \ll C(Y) i(S_X,L_Y)$.  It follows that
\begin{displaymath}
	i(S_X^u,L_X^u) \ll C(Y)/i(L_X,L_Y) .
\end{displaymath}
As $X \arrow \infty$ in $\cT_g$, we have $i(L_X,L_Y) \arrow \infty$ and hence
$i(S_X^u,L_X^u) \arrow 0$, which implies that the limiting laminations above satisfy $i(C,C') = 0$.
\qed

\bold{Proof of Theorem \ref{thm:pmlyes}.}
If $C$ is filling and uniquely ergodic, then the condition $i(C,C')$ implies that
$[C] = [C'] \mem \proj \ML_g$.
\qed

The following statement is elementary.
\begin{prop}
\label{prop:simplex}
For randomly chosen point $x$ in the simplex $\Delta \subset \reals_+^n$
defined by $\sum (x_i/A_i) \le 1$, we have $E(x_i) = A_i/(n+1)$.
\end{prop}

\bold{Proof of Theorem \ref{thm:pmlno}.}
Let $\sum_1^{3g-3} C_i \mem \ML_g(\zed)$ be the current associated to a maximal
system of disjoint simple closed $\alpha = \bigcup_1^{3g-3} \alpha_i \subset \Sigma_g$.

Extending $\alpha$ to a hexagon decomposition $(\alpha,\beta,H)$,
we obtain compatible Dehn coordinates $(m_i,t_i)$ on $\ML_g$ and 
Fenchel--Nielsen coordinates $(\ell_i,\tau_i)$ on $\cT_g$.

Choose weights $a_i>0$, and let $X_n$ be the point
in $\cT_g$ with Fenchel--Nielsen coordinates
\begin{displaymath}
	(\ell_i,\tau_i) = (1/(na_i),0)
\end{displaymath}
for $i=1,2,\ldots,3g-3$.  It is well--known that 
\begin{displaymath}
	[\lambda(X_n)] \arrow \left[\sum_1^{3g-3} C_i \right]
\end{displaymath}
as $n \arrow \infty$, since the width of the standard collar about $\alpha_i$ on $X_n$ is asymptotic to $2 \log n$, independent of $i$.

Now consider $\sigma(X_n)$.
By Theorem \ref{thm:Msharp}, for $n$ large the unit ball $B_X \subset \ML_g$ is very well
approximated in Dehn coordinates by the region where
\begin{displaymath}
	\sum |t_i| \ell_i + 2 m_i \Log(1/\ell_i) \le 1.
\end{displaymath}
Applying Proposition \ref{prop:simplex} to this simplex, we find that the coordinates of a random point $C =D(m_i,t_i) \mem B_X$ satisfy
\begin{displaymath}
	E(|t_i|) \asyto n a_i \AND E(m_i) \asyto 1/\log(n) .
\end{displaymath}
It follows from Theorem \ref{thm:i} that for any fixed multicurve $C'$ with Dehn coordinates $(m_i', t_i')$, we have
\begin{displaymath}
	E(i(C,C')) \asyto E\left(\sum m_i' |t_i|\right) \asyto n \sum a_i m_i' 
\end{displaymath}
as $n \arrow \infty$ (unless all $m_i'=0$, in which case $i(C,C') \arrow 0$).
Since $i(\sum a_i C_i,C') = \sum a_i m_i'$ as well, we have
\begin{displaymath}
	[\sigma(X_n)] \arrow \left[\sum a_i C_i\right] \;\;\text{in $\proj \cC_g$}, 
\end{displaymath}
completing the proof.
\qed

\bibliography{math}
\bibliographystyle{math}

\bigskip
{\sc Mathematics Department, Harvard University, 1 Oxford St, Cambridge, MA 02138}

\bigskip
{\sc Mathematics Department, University of Chicago,  5734 S University Ave, Chicago, IL 60637}

\end{document}